\documentclass{amsart}
\usepackage{amssymb,latexsym,graphics,amsfonts}

\swapnumbers
 \theoremstyle{plain}

 \theoremstyle{definition}
 
 \theoremstyle{remark}

 \numberwithin{equation}{subsection}

 \newcommand{\cal}[1]{\mathcal{#1}}

\begin{document}

\title{Derivations into n-th duals of ideals of Banach algebras}

\author{Madjid Eshaghi Gordji}
\address{Department of Mathematics,
University of Semnan, Semnan, Iran} \email{maj\_ess@Yahoo.com}

\author{ Reza Memarbashi}
\address{Department of Mathematics,
University of Semnan , Semnan, Iran}
\email{rmemarbashi@yahoo.com}



\dedicatory{}



\smallskip

\begin{abstract}
We introduce two notions of amenability for a Banach algebra $\cal
A$. Let $n\in \Bbb N$ and let $I$ be a closed two-sided ideal in
$\cal A$, $\cal A$ is $n-I-$weakly amenable if the first cohomology
group of $\cal A$ with coefficients in the n-th dual space $I^{(n)}$
is zero; i.e., $H^1({\cal A},I^{(n)})=\{0\}$. Further, $\cal A$ is
n-ideally amenable if $\cal A$ is $n-I-$weakly amenable for every
closed two-sided ideal $I$ in $\cal A$. We find some relationships
of $n-I-$ weak amenability and $m-J-$ weak amenability for some
different m and n or for different closed ideals $I$ and $J$ of
$\cal A$.
\end{abstract}

\maketitle


\section{Introduction}

Let $\cal A$  be a Banach algebra and $X$ a Banach $\cal A$-module,
that is $X$ is a Banach space and an $\cal A$-module such that the
module operations $(a,x)\longmapsto ax$ and $(a,x)\longmapsto xa$
from ${\cal A}\times X$ into $X$ are jointly continuous. Then $X^*$
is also a Banach $\cal A$-module if we define
\begin{center}
$\langle x,ax^*\rangle =\langle xa,x^* \rangle$\\
$\hspace{8cm}(a\in {\cal A}, ~x\in X,~x^*\in X^*)$\\
$\langle x,x^*a\rangle=\langle ax,x^*\rangle. $
\end{center}
In particular $I$ is a  Banach $\cal A$-module and $I^*$ is a dual
$\cal A$-module for every closed two-sided ideal $I$
of $\cal A$.\\

\noindent If $X$ a Banach $\cal A$-module then a derivation from
$\cal A$ to $X$ is a continuous linear operator $D$ with
$$D(ab)=a\cdot D(b)+D(a)\cdot b \hspace{2cm}(a,b\in {\cal A})~~~~(1.1).$$
If $x\in X$ and we define $D$ by
$$D(a)=a\cdot x-x\cdot a\hspace{1.5cm}(a\in {\cal A})$$
then $D$ is a derivation. Such derivations are called inner. A
Banach algebra $\cal A$  is amenable if every derivation from $\cal
A$ to every dual $\cal A$-module be inner; this definition was
introduced by B. E. Johnson in [Jo1] (see [Ru] and [He]). $\cal A$
is weakly amenable if, $H^1({\cal A},{\cal A}^*)=\{0\}$ (see [Jo3],
[D-Gh], [Gr1], [Gr2] and [Gr3]). Bade, Curtis and Dales [B-C-D] have
introduced the concept of weak amenability for commutative Banach
algebras. Let $n\in \Bbb N$. A Banach algebra $\cal A$ is called
n-weakly amenable if, $H^1({\cal A},{\cal A}^{(n)})=\{0\}.$ Dales,
Ghahramani and Gronbaek started the concept of n- weak amenability
of Banach algebras in [D-Gh-G]. A Banach algebra $\cal A$ is ideally
amenable if $H^1({\cal A},I^*)=\{0\}$ for every closed ideal $I$ of
$\cal A$ (see [G-Y]). In this paper we shall study $H^1({\cal
A},I^{(n)})$ for closed ideal $I$ of $\cal A$. Of course the
following definition describes
the main new property in this work.\\

\noindent {\large\bf Definition 1.1.} Let $\cal A$ be a Banach
algebra, $n\in \Bbb N$ and $I$ be a closed two-sided ideal in $\cal
A$. Then $\cal A$ is $n-I-$weakly amenable if $H^1({\cal
A},I^{(n)})=\{0\}$; $\cal A$
 is n-ideally amenable if $\cal A$ is $n-I-$weakly amenable   for every
 closed two-sided ideal $I$ in $\cal A$ and $\cal A$
 is permanently ideally amenable if $\cal A$ is $n-I-$weakly amenable   for every
 closed two-sided ideal $I$ in $\cal A$ and every $n\in \Bbb N.$\\

\noindent {\large\bf Example 1.} Let ${\cal  A}=\ell^1(\Bbb N)$. We
define the product on  ${\cal  A}$ by $f\cdot g=f(1)g \hspace{.2cm}
f,g\in  {\cal  A}$. ${\cal  A}$ is a Banach algebra with this
product and norm $\parallel\cdot\parallel _1$. Let $I$ be a closed
two-sided ideal of ${\cal A}$ it is easy to see that if $I\neq {\cal
A}$, then $I\subseteq \{f\in {\cal  A};~f(1)=0\}$. Then The right
module action of $\cal A$ on $I$ is trivial, therefore the right
module action of $\cal A $ on $I^{(2k)}$ is trivial for every $k\in
\Bbb N$. On the other hand $\cal A$ has left identity, then by
proposition 1.5 of [Jo1] we have  $H^1({\cal A},I^{(2k+1)})=\{0\}
(k\geq 0)$. If $I=\cal A$ then by Assertion 2 of [Zh], $H^1({\cal
A},I^{(2k+1)})=\{0\} (k\geq 0)$. Thus for every $k\geq 0$, $\cal A$
is 2k+1 ideally amenable. It is well known that $\cal A$ is not
2-weakly( ideally) amenable (see [Zh]).\\

Let $X$ be a Banach $\cal A-$module, we can extend the actions of
${\cal A}$ on $X$ to  actions of ${\cal A}^{**}$  on $X^{**}$ via
\[a''.x'' =w^{*}\textup{-}\lim_{i} \lim_{j} a_{i} \, x_{j}
\]
and
\[ x''.a'' =w^{*}\textup{-}\lim_{j} \lim_{i} x_{j} \, a_{i},\]
where $a'' =w^{*}\textup{-}\lim_{i} a_{i}$, \  $x''
=w^{*}\textup{-}\lim_{j} x_{j}$.

\noindent{\large\bf Definition.1.2.} Let $\cal A$ be a Banach
algebra and $X$ be a Banach $\cal A-$module. We define the
topological center of the right module action of $\cal A$ on $X$ as
follows
 $$Z_{\cal A}(X^{**}):=\{x''\in X^{**}~:~
 \hbox{the mapping } a''\mapsto
x''.a'':{\cal A}^{**}\rightarrow X^{**} \hbox{ is~$weak^*-weak^*$ ~
continuous }\}$$

The right module action of $\cal A$ on $X$ is Arens regular if and
only if $Z_{\cal A}(X^{**})=X^{**}$ (see [A] and [E-F]). For a
Banach algebra $\cal A$ the set $Z_{\cal A}(\cal A^{**})$ is the
topological center of $\cal A^{**}$ with the first Arens
product.\\
Let $\cal A$ be a Banach algebra and $X$ be a Banach $\cal
A-$module. Let $P:X''''\rightarrow X''$ be the adjoint of the
inclusion map $i:X'\rightarrow X'''.$ Then we have\\

\noindent {\large\bf Theorem 1.3.} Let $\cal A$ be a Banach algebra
and $X$ be a Banach $\cal A-$module. Suppose that $Z_{\cal
A}(X'')=X''$. Then the following assertions holds
\begin{quote}
 (i)
$P:X''''\rightarrow X''$ is an ${\cal A}^{**}-$module
morphism.\\
(ii) Let $D:\cal A \rightarrow X''$ be a derivation then there
exists a derivation $\tilde {D}:{\cal A}^{**} \rightarrow X''$ which
$\tilde {D}$ is extension of $D$.
\end{quote}
\noindent {\large\bf Proof.} (i) conclude from Proposition 1.8 of
[D-Gh-G]. For (ii), we know that $D'':{\cal A}^{**} \rightarrow
X''''$ the second adjoint of $D$ is a derivation (see for example
Proposition 1.7 of [D-Gh-G]). By (i) $P\circ D''$ is a derivation
from
${\cal A}^{**}$ into $X''$. \\

\noindent {\large\bf Corollary 1.4.} Let $\cal A$ be Arens regular
Banach algebra. If for every ideal $I$ of  ${\cal A}^{**}$
 $H^1({\cal A^{**}},I^{**})=\{0\}$, then $\cal
A$ is 2-ideally amenable.\\

 For convenience we will write $x\mapsto J(x)$ for
the canonical embedding of a Banach space into its second dual. We
fined some relations between m and n ideal amenability of a Banach
algebra.\\

\noindent {\large\bf Theorem 1.5.} Let ${\cal A}$ be a Banach
algebra and $I$ be a closed ideal of ${\cal A}$. For any $n\in \Bbb
N$, if $\cal A$ be $n+2-I-$weakly amenable then $\cal A$ is
$n-I-$weakly amenable.\\
 \noindent {\large\bf Proof} Let $D:{\cal A}\rightarrow I^{(n)}$ be a
 derivation. Since $J:I^{(n)}\rightarrow I^{(n+2)}$ is an $\cal A-$
 module homomorphism, then $J\circ D:{\cal A}\rightarrow I^{(n+2)}$ is a
 derivation. Then there exists $F\in I^{(n+2)}$ such that $J\circ D=\delta_F.$
Let $P:I^{(n+2)}\rightarrow I^{(n)}$ be the above projective, then
for every $a\in \cal A$ we have $D(a)=P\circ J\circ D(a)=a\cdot
P(F)-P(F)\cdot a$.
Therefore $D=\delta_{P(F)}.$\\

\noindent{\large\bf Corollary 1.6.} Let $\cal A$ be a Banach
algebra, and let $n\in \Bbb N.$ If $\cal A$ be $n+2-$ideally
amenable then $\cal A$ is $ n-$ideally amenable.
\\

 \noindent {\large\bf Theorem 1.7.} Let ${\cal A}$ be a Banach
algebra and $I$ be a closed two sided ideal of ${\cal A}$ with a
bounded approximate identity. If ${\cal A}$ is n-ideally amenable
(or permanently ideally amenable) then $I$ is
n-ideally amenable (or permanently ideally amenable).\\
\noindent {\large\bf Proof.} Since $I$ has bounded approximate
identity, then by Cohen factorization Theorem for every closed ideal
$J$ of $I$, we have $JI=IJ=J.$ Then $J$ is an ideal of ${\cal A}$.
Let $ D : I \longrightarrow J^{(n)}$ be a derivation. By [Ru,
Proposition 2.1.6], $D$ can extends to a derivation ${\bf \tilde D}
: {\cal A} \longrightarrow J^{(n)}$. So there is a $m \in J^{(n)}$
such that ${\bf\tilde D}= \delta_{m}$. Then $D(i)={\bf\tilde
D}(i)=\delta_m$ for each $i \in I$. Thus $D$ is
inner.\\

\noindent {\large \bf Corollary 1.8.} Let ${\cal A} $ be a Banach
algebra  with a bounded approximate identity and let ${\cal M(A)}$
the multiplier algebra of  $ {\cal A}$ be n-ideally amenable (or
permanently ideally amenable) then ${\cal A}$ is n-ideally
amenable (or permanently ideally).\\

\noindent {\large\bf Theorem 1.9.} Let ${\cal A}$ be a Banach
algebra and $I$ be a closed two sided ideal of ${\cal A}$ with unite
element. Let $J$ be a closed ideal of $I$, and $n\in \Bbb N$. Then
$J$ is a closed ideal of $\cal A$ also ${\cal A}$ is n-J-weakly
amenable if and only if $I$ is n-J-weakly amenable.\\
\noindent {\large\bf Proof.} Since $I$ has unite element, then
$J^{(n)}$ is essential. Then every derivation  $ D : I
\longrightarrow J^{(n)}$ has a unique extension $ \tilde{D} : {\cal
A} \longrightarrow J^{(n)}$ such that $\tilde{D}$ is a ( bounded)
derivation [D-Gh-G]. Thus we can prove the necessary and
sufficient conditions easily.\\

\noindent {\large\bf Theorem 1.10.} Let ${\cal A}$ be a Banach
algebra and let $n\in \Bbb N$. Let $I$ be a closed ideal of ${\cal
A}$ which $Z_{\cal A}(I^{(2n)})=I^{(2n)}.$ Suppose that $H^1({\cal
A^{**}},I^{(2n)})=\{0\}$, then $\cal A$ is 2n-2-I-weakly amenable.\\
\noindent {\large\bf Proof.} Let $ D : \cal A \longrightarrow
I^{(2n-2)}$ be a derivation then by Theorem 1.3, there exists an
extension $\tilde{D} : \cal A^{**} \longrightarrow I^{(2n)}$ such
that $\tilde{D}$ is a (bounded) derivation. Then $\tilde{D}$ is
inner so $D$ is inner.\\

\noindent {\large\bf Theorem 1.11.} Let ${\cal A}$ be a Banach
algebra with a left bounded approximate identity. Let $I$ be a
closed ideal of ${\cal A}$ and let $\cal A$ is ideal in ${\cal
A}^{**}$. If $I$ is left strongly irregular (i.e.$Z_t(I^{**})=I$)
and $\cal A$ is $I-$weakly amenable then $\cal A$ is $3-I-$weakly
amenable.\\
\noindent {\large\bf Proof.} First we have the following $\cal
A-$module direct sum decomposition
$$I^{***}=\hat {I^*}\oplus {\hat{I}}^{\perp}$$ then we have
$$H^1{(\cal A,I^{***})}= H^{1}{(\cal A,\hat {I^*})}+H^1{(\cal A,\hat
{I}^{\perp})}.$$ Thus it is enough  to show that $H^1{(\cal A,\hat
{I}^{\perp})}=\{0\}.$ Let $a\in \cal A$ and let $i''\in I^{**}$ and
$\pi:I\rightarrow \cal A$ be the inclusion map. Then by Lemma 3.3 of
[Gh-L],
$$i''a=\pi''(i'')\hat {a}\in \pi''(I^{**})\cap \hat {\cal A}=\hat
{I}$$ Then the right module action of $\cal A$ on $\hat {I}^{\perp}$
is trivial. Let now $D:\cal A\rightarrow \hat {I}^{\perp}$ be a
derivation. Suppose $(e_\alpha)$ be a left bounded approximate
identity for $\cal A.$ Since $\hat {I}^{\perp}$ is a $weak^*-$closed
sub space of $I^{***},$ we take $weak^*-\lim_\alpha D(e_\alpha)=F\in
\hat {I}^{\perp}.$ Then for every $a\in \cal A$ we have
$$D(a)=\lim_\alpha D(e_\alpha a)=Fa=Fa-aF=\delta_F(a).$$

\section{ The unitization of a Banach algebra}

Let ${\cal A}^{\#} $ be the unitization of ${\cal A}$. We know that
${\cal A}$ is amenable if and only if ${\cal A}^{\#}$ is amenable.
If ${\cal A}$ is weakly amenable then ${\cal A}^{\#}$ is weakly
amenable, [D-Gh-G]. and the weak amenability of ${\cal A}^{\#} $ dos
not implies the weak amenability of ${\cal A}$ [Jo-Wh]. Also the
first author  and T.Yazdanpanah shown that ${\cal A}$ is ideally
amenable if and only if ${\cal A}^{\#}$ is ideally amenable. In the
following we will prove the same problem for
n-ideal amenability.\\

\noindent {\large\bf Proposition 2.1.} Let ${\cal A}$ be a Banach
algebra, and let $n\in \Bbb N.$ Then the following assertions hold:
\begin{quote}
(i) If ${\cal A}^{\#}$ be n-ideally amenable then ${\cal A}$ is
n-ideally amenable.
\\
(ii) If ${\cal A}$ be $2n-1-$ideally amenable then ${\cal A}^{\#}$
is $2n-1-$ideally amenable.
\\
(iii) If ${\cal A}$ be commutative and n-ideally amenable, then
${\cal A}^{\#}$ is n-ideally amenable.
 \end{quote}
\noindent {\large\bf Proof. } (i) Let ${\cal A}^{\#}$ be n-ideally
amenable, and let $I$ be a closed ideal of ${\cal A}$ and $ D :
{\cal A} \longrightarrow I^{(n)}$ be a derivation. It is easy to see
that $I$ is an ideal of ${\cal A}^{\#}$. We define ${\bf\tilde D} :
{\cal A}^{\#} \longrightarrow I^{(n)}$ by  $ {\bf\tilde D}
(a+\alpha) = D(a)$, $ (a \in {\cal A} ,~ \alpha \in \Bbb C ).$ Then
$\bf\tilde D$ is a derivation. Since ${\cal A}^{\#}$ is n-ideally
amenable then $ \bf\tilde{D} $ is inner
and hence $D$ is inner.\\
For (ii), let $ {\cal A}$ be $2n-1-$ideally amenable and $I$ be a
closed ideal of ${\cal A}^{\#}$. First we know that ${\cal A}$ is
$2n-1-$weakly amenable. Then by proposition 1.4 of [D-Gh-G], ${\cal
A}^{\#}$ is $2n-1-$weakly amenable. Thus ${\cal A}^{\#}$ is
$2n-1-I-$weakly amenable whenever $I={\cal A}^{\#}.$ Let now $I\neq
{\cal A}^{\#},$ then $1 \notin I$ and $I$ is an ideal of ${\cal A}$.
Let $ D : {\cal A}^{\#} \longrightarrow I^{(n)}$  be a derivation,
then $ D(1)=0$ and $D$ drops to a derivation from ${\cal A}$ into
$I^{(n)}$ and $D$ is
inner.\\
The proof of (iii) is similar (ii). \\

\noindent {\large\bf Proposition 2.2.} Let ${\cal A}=L^1(G)$ for
locally compact group discrete  $G$ and let $M$ be the augmentation
ideal of $\cal A$, then $\cal A=M^{\#}$ therefore
\begin{quote}
(i) If ${\cal A}$ be n-ideally amenable then $M$ is n-ideally
amenable for every $n\in \Bbb N.$
\\
(ii) If $M$ be $2n-1-$ideally amenable then ${\cal A}$ is
$2n-1-$ideally amenable for every $n\in \Bbb N$.
\end{quote}

\section{ Commutative Banach algebras.}

We know that a commutative Banach algebra is weakly amenable if and
only if every derivation from $\cal A$ into a commutative Banach
$\cal A$-module is zero (Theorem 1.5 of [B-C-D]). Then we conclude the following.\\

\noindent {\large\bf Theorem 3.1.} Let $\cal A$ be a commutative
Banach algebra, then the following assertions are equivalent.
\begin{quote}
(i) $\cal A$ is weakly amenable.\\
(ii) $\cal A$ is 2k+1-weakly amenable for some $k\in \Bbb N \cup
\{0\}$.\\
(iii) $\cal A$ is ideally amenable.\\
(iv)  $\cal A$ is 2k+1-ideally amenable for some $k\in \Bbb N \cup
\{0\}$.\\
(v) $\cal A$ is permanently ideally amenable.\\
\end{quote}

\noindent {\large\bf Theorem 3.2.} Let $\cal A$ be a commutative
Banach algebra and let $n\in \Bbb N$, then the following assertions
are equivalent.
\begin{quote}
(i) $\cal A$ is 2n-weakly amenable.\\
(ii) $\cal A$ is 2n-ideally amenable.
\end{quote}
{\large\bf Proof.}  (ii) $\Rightarrow$ (i) is obviously. For (i)
$\Rightarrow$ (ii), let $\cal A$ be 2n-weakly amenable and let $I$
be a closed two sided ideal of $\cal A$. We denote
$\pi:I\longrightarrow \cal A$ the natural inclusion map, then
$\pi^{(2n)}:I{(2n)}\longrightarrow {\cal A}^{(2n)}$ the 2n-th
adjoint of $\pi$  is $\cal A-$module morphism. Let $ D : {\cal A}
\longrightarrow I^{(2n)}$ be a derivation, then $\pi^{(2n)}oD :
{\cal A} \longrightarrow {\cal A}^{(2n)}$ is a derivation. Since
${\cal A}^{(2n)}$ is commutative $\cal A-$module and $H^1({\cal
A},{\cal A}^{(2n)})=\{0\},$ then $\pi^{(2n)}oD=0.$ Therefore $D=0.$
\\

\noindent {\large\bf Corollary 3.3.} Let $\cal A$ be a commutative
Banach algebra which is Arens regular, and suppose that ${\cal
A}^{**}$ is semisimple. Then $\cal
A$ is 2-ideally amenable.\\
\noindent {\large\bf Proof.} By Corollary 1.11 of [D-Gh-G], $\cal A$
is 2-weakly amenable. Then by Theorem 3.2 $\cal A$ is 2-ideally
amenable.\\

\noindent {\large\bf Corollary 3.4.} Let $\cal A$ be a commutative
Banach algebra such that ${\cal A}^{(2n)}$ is  Arens regular, and
$H^1({\cal A}^{(2n+2)}, {\cal A}^{(2n+2)})=\{0\}$ for each $n\in
\Bbb N$. Then $\cal A$ is 2n-ideally amenable for each $n\in\Bbb
N$.\\
\noindent {\large\bf Proof.} By Corollary  1.12 of [D-Gh-G]  $\cal
A$ is 2n-weakly amenable for each $n\in\Bbb N$. Then  by Theorem 3.2
$\cal A$ is 2n-ideally amenable for each $n\in\Bbb N$.\\

\noindent {\large\bf Corollary 3.5.}  Every uniform Banach algebra
is 2n-ideally  amenable for each $n\in\Bbb N$.\\
{\large\bf Proof.} By applying   Theorem 3.2 above and Theorem 3.1
of
[D-Gh-G] the proof if easy.\\

Let $\Bbb D=\{z\in\Bbb C~:~|z|<1\}$ the open unit disc and let
$A(\overline {\Bbb D})$  be the disc algebra. It follows from
Corollary 3.5 above and page 35 of [D-Gh-G] that $A(\overline {\Bbb
D})$ is a 2-ideally amenable Banach function algebra which is not
ideally amenable.

\section{${\cal C}^*-$algebras.}

 It is well known that every  ${\cal
C}^*-$algebra is ideally amenable [G-Y; Corollary 2.2.] also a
${\cal C}^*-$algebra is amenable if and only if it is nuclear
([Ha]). \\
As in [D-Gh-G] we have the following.\\

\noindent {\large\bf Theorem 4.1. } Every $C^*$-algebra is
Permanently weakly amenable.\\
\vspace{-.38cm}

We can not show that every $C^*-$algebra is permanently ideally
amenable but we prove the following Theorem for $C^*-$algebras.\\

\noindent {\large\bf Theorem 4.2. } Let $n=2k+1~~~ (k\in {\Bbb
N}\bigcup \{0\}).$ Then every $C^*$-algebra is
n-ideally amenable.\\
\noindent {\large\bf Proof.} Let ${\cal A}$ be a $C^*$-algebra and
let $I$ be a closed ideal of ${\cal A}$. Since $\cal A$ is ideally
amenable then $H^1({\cal A},I^{*})=\{0\}$. Now we show that if every
$C^*$-algebra be $n=2k+1-$ ideally amenable then every $C^*$-algebra
is $2k+3-$ ideally amenable. Let $k\geq 0,$ and let $D:{\cal
A}\rightarrow
 {I^{(n+2)}}$ be a
 derivation. Then
 we show that $D''$ is a derivation. Let $a'',b''\in
{\cal A}^{**}$ then there are nets $(a_\alpha)$ and $(b_\beta)$ in
${\cal A}$ such that converge respectively to $a''$ and $b''$ in the
$weak^*$- topology of ${\cal A}^{**}$. Then
\begin{eqnarray*}
D''(a''b'')&=& weak^* lim_\alpha lim_\beta D(a_\alpha b_\beta) \\
&=&weak^* lim_\alpha lim_\beta D(a_\alpha)b_\beta+weak^* lim_\alpha
lim_\beta a_\alpha D(b_\beta) \\
&=& D''(a'').b''+ lim_\alpha  a_\alpha.
D''(b'')~\hspace{2.5cm}~~\hfill~~~~~~(4.1).
\end{eqnarray*}
Let $x''\in {I^{(n+3)}}$ and let $\pi:I\rightarrow \cal A$ be the
inclusion map and let $i:\cal A\rightarrow {\cal A}^{**}$ be the
natural embedding. The maps $i'', i^{(4)},...,i^{(n+3)} $ are
$weak^*$-$weak^*-$continuous. Then $~~weak^*-\lim_\alpha
i^{(n+3)}(a_\alpha)=i^{(n+3)}(a'').~~~~$ On the other hand ${\cal A}^{(n+3)}$ is a $C^*-$algebra, then it is Arens regular, thus \\
$$\lim_\alpha x''a_\alpha=\lim_\alpha
\pi^{(n+3)}(x'')i^{(n+3)}(a_\alpha)=\pi^{(n+3)}(x'')i^{(n+3)}(a'')=x''a''$$
Since every bounded linear map from every $C^*-$algebra is weakly
compact, then $D$ is weakly compact, therefore  $D''(b'')\in \hat
{I^{(n+2)}}$, then for each $x''\in {I^{(n+3)}}$, we have
\begin{eqnarray*}
lim_\alpha \langle a_\alpha.D''(b''),x''\rangle&=&lim_\alpha \langle
x'' a_\alpha,D''(b'')
\rangle \\
&=&\langle x''a'',D''(b'')\rangle \\
&=&\langle a''.D''(b''),x''\rangle.
\end{eqnarray*}
Then
 $$a''.D''(b'')=lim_\alpha  a_\alpha. D''(b'').$$
and by 4.1,  $D''$ is a derivation. Since $D$ is weakly compact,
then $D''({\cal A}^{**})\subseteq \hat {I^{(n+2)}}.$ We can suppose
that $D''$ is a derivation from $\cal A^{**}$ into $I^{(n+2)}$.
Similarly $D^4:=(D'')^{''}$ is a derivation from ${\cal A}^{(4)}$
into $hat {I^{(n+2)}}$. We can suppose that $D^{(2k+2)}$ (the
$2k+1-$th conjugate of $D$) is a derivation from ${\cal A}^{(2k+2)}$
into $\hat {I^{(n+2)}}$. On the other hand $I^{(n+1)}$ is a closed
ideal of ${\cal A}^{(2k+2)}$ and ${\cal A}^{(2k+2)}$ is ideally
amenable (since it is a $C^*-$algebra). Thus $D^{(2k+2)}$ is inner.
Since $D^{(2k+2)}$ is extension of $D$, then
it is easy to see that $D$ is inner.\\

\section{ codimension one ideals}
Let $\cal A$  be a Banach algebra and let $I$ be a closed ideal of
$\cal A$ by codimension one. We fined the relationship between
n-weak amenability of $I$ and $n-I-$weak amenability of $\cal A.$
As [7;Theorem. 2.3] we have the following.\\

\paragraph{\large\bf Theorem 5.1.} Let $\cal A$ be a Banach algebra with bounded approximate
identity and let $I$ be a codimension one closed two sided ideal of
$\cal A$. Then $H^1({\cal A},X^*)\cong H^1({\cal I},X^*)$ for every
neo unital Banach $\cal A-$
module $X$.\\

Let $G$ be a discrete group, and let $I_0$ be a codimension 1 closed
two sided ideal of $l^1(G)$. Then
$l^1(G)$ is $n-I_0-$weakly amenable for every $n\in \Bbb N$.\\

\paragraph{\large\bf Corollary 5.2.} Let $\cal A$ be a $C^*-$algebra
and let $I$ be a codimension one closed two sided ideal of $\cal A$.
Then for every $n\in \Bbb N$, $\cal A$ is $n-I-$weakly amenable.
\paragraph{\large\bf Proof.} Let $n=2k+1$ then
by Theorem 4.2, $H^1({\cal A},I^{(n)})=\{0\}$. Let $n=2k$, then we
have ${\cal A} I^{(n-1)}=II^{(n-1)}=I^{(n-1)}$ and $I^{(n-1)} {\cal
A}=I^{(n-1)}I=I^{(n-1)}.$ Then by Theorem 1.5 $\cal A$ is
$n-I-$weakly amenable.\\

\paragraph{\large\bf Lemma 5.3.} Let $\cal A$ be a Banach algebra with bounded
approximate identity and let $I$ be a closed two sided ideal of
$\cal A$. For every $n\geq 0$, the following assertions hold.\\
(i) If $\pi_r:I^{(n)}\times{\cal A}\longrightarrow I^{(n)}$ be Arens
regular,
then $I^{(n+1)}$ factors on the right.\\
(ii)If $\pi_l:{\cal A}\times I^{(n)}\longrightarrow I^{(n)}$ be
Arens regular, then $I^{(n+1)}$ factors on the left.
\paragraph{\large\bf Proof.} Let $(e_\alpha)$ be a bounded approximate
identity for $\cal A$ with cluster point $E$. Suppose that $\pi_r$
be Arens regular and let $i^{n+2}\in I^{(n+2)}$ be the cluster point
of $(i^n_\beta)$ ($(i^n_\beta)$ is a net in $I^{(n)}$). Then we have
\begin{eqnarray*}
i^{n+2}&=&\lim_\beta i^{n}_\beta \\
&=&\lim_\beta\lim_\alpha i^n_\beta e_\alpha  \\
&=&\lim_\alpha\lim_\beta i^n_\beta e_\alpha \\
&= &\lim_\alpha i^{(n+2)} e_\alpha =i^{(n+2)}E.
\end{eqnarray*}
Then for $i^{n+1}\in I^{(n+1)}$ we have
\begin{eqnarray*}
\lim_\alpha\langle e_\alpha
i^{n+1},i^{n+2}\rangle&=&\lim_\alpha\langle i^{n+1},
i^{n+2}e_\alpha\rangle   \\
&=&\langle i^{n+1}, i^{n+2}E \rangle  \\
&=&\langle i^{n+1},i^{n+2}\rangle.
\end{eqnarray*}
Thus $e_\alpha i^{n+1} \rightarrow i^{n+1}$ weakly. Since $e_\alpha
i^{n+1} \in {\cal A}I^{(n+1)}$ for every $\alpha$, by Cohen-Hewit
factorization theorem we know that ${\cal A}I^{(n+1)}$ is closed in
$I^{(n+1)}$, then $i^{n+1} \in {\cal A}I^{(n+1)}$. Thus the proof of
(i) is complete. For (ii), let $i^{n+2}\in I^{(n+2)}$ be the cluster
point of $(i^{n}_\beta)$, for each $\beta$ we know that
$Ei^{n}_\beta=i^{n}_\beta$ since $\pi_l$ is Arens regular, then
$e_\alpha i^{n+2} \rightarrow Ei^{n+2}=i^{n+2}$ by $weak^*$ topology
of $I^{(n+2)}$. Then for every $i^{n+1} \in I^{(n+1)}$,
\begin{eqnarray*}
\lim_\alpha\langle i^{n+1}e_\alpha
,i^{n+2}\rangle&=&\lim_\alpha\langle i^{n+1},
e_\alpha i^{n+2}\rangle   \\
&=&\langle i^{n+1}, Ei^{n+2} \rangle  \\
&=&\langle i^{n+1},i^{n+2}\rangle.
\end{eqnarray*}
Therefore $i^{n+1}e_\alpha \rightarrow i^{n+1}$ weakly. Again by
Cohen-Hewit
factorization theorem we conclude that $I^{(n+1)}$ factors on the left.\\
By Theorem 5.1 and Lemma 5.3, we have the follows\\

\paragraph{\large\bf Proposition 5.4.} Let $\cal A$ be a Banach algebra with bounded
approximate identity and let $I$ be a codimension one closed two
sided ideal of $\cal A$. For every $n\geq 1$, if the module actions
of $\cal A$ on $I^{(n-1)}$ are Arens regular, then  $I$ is
$n-$weakly amenable if and only if $\cal A$
is $n-I-$weakly amenable.\\

\section{ Open Problems}
\begin{itemize}

\item We do not know whether or not 2-ideal amenability implies
4-ideal amenability for arbitrary Banach algebra.

\item we do not know whether or not 1-ideal amenability implies 3-ideal
amenability for arbitrary(non commutative) Banach algebra.

\item We know that every $C^*-$algebra is 2k+1 ideally amenable for every
$k\in \Bbb N \bigcup \{0\}.$ But we do not know whether or not every
$C^*-$algebra is permanently ideally amenable.

\item The group algebras $L^1(G)$ are n-weakly amenable for each odd n
[D-Gh-G], but we do not know for which $G$
 and which $n\in \Bbb N$, the algebra $L^1(G)$ is n-ideally
 amenable. We know that $l^1(\Bbb F_2)$ the group algebra of free group $\Bbb
F_2$ is permanently amenable [Jo2], but we do not know whether or
not $l^1(\Bbb F_2)$ is permanently ideally amenable.

\item We know that sometimes weak amenability of ${\cal A}^{**}$
implies the weak amenability of $\cal A$ (see [GH-L-W], [Gh-La] and
[D-G-V]). But we do not know when the ideal amenability of ${\cal
A}^{**}$ implies the ideal amenability of $\cal A?$
\end{itemize}





\end{document}